\documentclass[12pt]{article}
\usepackage{graphicx}
\usepackage{psfrag}
\usepackage{color}
\usepackage{amssymb,latexsym}

\def\L{\left}
\def\R{\right}
\def\Ref#1{(\ref{#1})}
\newcommand{\E}{\mathbf E}

\newcommand{\RR}{\mathbf R}

\renewcommand{\Pr}{\mathbf P}

\newtheorem{lemma}{Lemma}
\newtheorem{proposition}{Proposition}
\newtheorem{corollary}{Corollary}
\newtheorem{theorem}{Theorem}

\def\AA{{\cal A}(d,\alpha)}

\begin{document}
\title{Comparison results and steady states for the \\
Fujita equation with fractional Laplacian\\[1em]
R\'{e}sultats de comparaison et solutions stationnaires 
de l'\'{e}quation de Fujita avec Laplacien fractionnel
}

\author{
{\sc Matthias Birkner\footnote{
FB Mathematik, J.W. Goethe Universit\"{a}t,
D-60054 Frankfurt am Main, Germany.} ,   
Jos\'{e} Alfredo L\'{o}pez-Mimbela\footnote{Centro de Investigaci\'{o}n en Matem\'{a}ticas, Apartado Postal 402, 36000 Guanajuato, Mexico.\newline\vskip-.275cm\indent\hskip.3cm
{\sc Key words and phrases:}  
Blow-up and extinction of solutions of semilinear PDEs, \hfill compar-\newline\indent\hskip.3cm ison,
Feynman-Kac representation, symmetry of solutions,
symmetric stable processes, method\newline\indent\hskip.3cm of moving planes.\newline\vskip-.3cm\indent\hskip.3cm
{\sc AMS subject classification:}
60H30, 35K57, 35B05, 60G52.
 }
} \\ and {\sc Anton Wakolbinger}\setcounter{footnote}{0}\footnotemark}
\date{ }

\maketitle
\begin{abstract}
We study a semilinear PDE generalizing the Fujita equation whose 
evolution operator is the sum of a fractional power of the 
Laplacian and a convex non-linearity. 
Using the Feynman-Kac representation we prove criteria for 
asymptotic extinction versus finite time blow up 
of positive solutions based on comparison with global solutions. 
For a critical power non-linearity we obtain 
a two-parameter family of radially symmetric stationary solutions. 

By extending the method of moving planes to fractional powers of 
the Laplacian we prove that all positive steady states of the 
corresponding equation in a finite ball are radially 
symmetric. 
\begin{center} \bf R\'{e}sum\'{e}\end{center}
Nous \'{e}tudions une \'{e}quation de r\'{e}action-diffusion 
semilin\'{e}aire
 (g\'{e}n\'{e}ralisant l'\'{e}qua\-tion de Fujita), dont l'op\'{e}rateur 
d'\'{e}volution
 est la somme d'une puissance fractionnelle du Laplacien et d'une
 non-lin\'{e}arit\'{e} convexe. 
A l'aide de la repr\'{e}sentation de Feynman-Kac
 nous exhibons des crit\`{e}res entrainant l'extinction asymptotique,
 respectivement l'explosion en temps fini, de solutions positives. 
Ces crit\`{e}res s'obtiennent en comparant avec des solutions globales.
 Pour une certaine puissance critique de la non-lin\'{e}arit\'{e} 
nous obtenons 
une famille param\'{e}tris\'{e}e de solutions stationnaires \`{a}  
sym\'{e}trie radiale.

 Par extension de la m\'{e}thode de d\'{e}placement d'hyperplans 
\`{a} des puissances
 fractionnelles du Laplacien, nous prouvons que toute solution
 positive stationnaire de l'\'{e}quation correspondante dans une boule
 finie comporte une  sym\'{e}trie radiale. 
\end{abstract}

\section{Introduction}

We consider the ``generalized Fujita equation''
\begin{eqnarray}
\label{IVP} \frac{\partial}{\partial t} u(t,x) & = & Lu(t,x) +
G(u(t,x)),\ \ t\ge0,\ x\in\RR^d, \\ \nonumber u(0,x) & = &
\varphi(x)\ge0,
\end{eqnarray}
where $L=\Delta_{\alpha}$ is the fractional power
$-(-\Delta)^{\alpha/2}$ of the Laplacian, $0<\alpha\le 2$, and
$G:\RR_+\to\RR_+$ is a convex function satisfying conditions
\Ref{G1} and \Ref{G2} below. 
Solutions will be understood in the mild sense 
(see e.g. \cite{Pazy}) so that \Ref{IVP} makes sense for 
any non-negative bounded measurable function $\varphi$ on $\RR^d$.

A well-known fact is 
that for any non-trivial initial value $\varphi$ there exists a number
$T_{\varphi} \in (0, \infty]$ such that \Ref{IVP} has a
unique solution $u$ on $\RR^d\times[0,T_{\varphi})$ which is
bounded on $\RR^d\times[0,T]$ for any $0<T<T_{\varphi}$, and if
$T_{\varphi}<\infty$, then
$\|u(\cdot,t)\|_{L^{\infty}(\RR^d)}\to\infty$ as $t\uparrow
T_{\varphi}$.

When $T_{\varphi}=\infty$ we say that $u$ is a global solution,
and when $T_{\varphi}<\infty$ we say that $u$ blows up in finite
time or that $u$ is non-global.

The study of blow up properties of \Ref{IVP} goes back to the
fundamental work of Fujita \cite{Fu}, who studied  Eq. \Ref{IVP}
with $\alpha=2$ and $G(z)=z^{1+\beta}$, $\beta>0$. The
investigation of \Ref{IVP} with a general $\alpha$ was initiated
by Sugitani \cite{Su}, who showed that if $d\le\alpha/\beta$, then for any
non-vanishing initial condition the solution blows up in finite time. 
Using a Feynman-Kac representation for the solutions of 
semilinear problems of the form \Ref{IVP}, 
this conclusion was re-derived in \cite{BLMW} and the 
corresponding behaviors of equations with time-dependent nonlinearities and 
of various systems of semilinear pde's were studied. 

It is known (e.g. \cite{LMW}, \cite{NS}) that in supercritical dimensions 
$d>\alpha/\beta$, Eq. \Ref{IVP} admits global as well as 
non-global positive solutions, 
depending on the ``size'' of the initial condition. 
For short we address this parameter constellation as the {\em global regime.}

In the first part of the present note we prove two comparison
criteria in the global regime:

(i) Assume the initial value $\varphi\ge 0$ leads to a globally bounded
solution. Then any initial value $\psi$ with  $0 \le \psi \le
(1-\varepsilon)\varphi$, $\varepsilon > 0$,
gives rise to a solution converging to zero.

(ii) Assume the initial value $\varphi\ge 0$ leads to a
solution which is uniformly bounded away from $0$ for all $t>0$ and
all $x$ in some ball $\subset \RR^d$. Then any initial
value $\psi$ with  $\psi \ge
(1+\varepsilon)\varphi$, $\varepsilon > 0$,
gives rise to a solution which blows up in finite time.

The essential tool in proving (i) and (ii) is the probabilistic
representation of the solution of (\ref{IVP}) provided by the
Feynman-Kac formula, that was obtained in \cite{BLMW} (see (\ref{FK4u}) below).

Natural candidates for the comparison in (i) and (ii) are
(time-)stationary solutions of (\ref{IVP}), i.e. solutions of the
``elliptic'' equation
\begin{equation}\label{ell}
\Delta_{\alpha}u(x) + G(u(x))=0, \quad x \in \RR^d.
\end{equation}
In the case $\alpha =2$ and $G(z)=z^{1+\beta}$, it is known
that (see \cite{GS}, \cite{CL}, \cite{GNN}, \cite{Wang})

- for $d > 2$, $1+\beta  < (d+2)/(d-2)$,
apart from $u \equiv 0$, no bounded non-negative
 solution of (\ref{ell}) exists 

- for $d > 2$,  $\beta = (d+2)/(d-2)-1$, all bounded
stationary solutions of (\ref{ell}) are given by the
family
\begin{equation}\label{family}
u_{c,A}(x) = \frac{A(d(d-2))^{(d-2)/2}}{\L(d(d-2) + \L(A^{2/(d-2)}\|x-c\|\R)^2 \R)^{(d-2)/2}}, \ \ \ c,x\in \RR^d, \ \ A\in \RR_+.
\end{equation}
(Note that the two parameters of the family are the symmetry center
$c$ of $u$ and its value $A$ at $c$.)

- for $d > 2$, $\beta > (d+2)/(d-2)-1$, there exists a one-parameter
family $u_{A}, A>0$ of solutions of (\ref{ell}) with the
properties: $u_{A}(x)$ is symmetric around $x=0$,
$$u_{A}(0)=A, \quad  \|x\|^{2/\beta}u_{A}(\|x\|)\to
K(d,\beta)
\mbox { as }\|x\| \to \infty$$ where $K(d,\beta)$ is a constant not depending on
$A$.

In the case $\alpha < 2$, much less is known. For $d > \alpha$,
$\beta = (d+\alpha)/(d-\alpha)-1$, we specify in Proposition
\ref{steadystates} a
two-parameter family $u_{c,A}$, $c\in\RR^d$, $0<A<\infty,$ of radially symmetric
solutions of
\begin{equation}\label{ellbeta}
\Delta_{\alpha}u(x) + u^{1+\beta}(x)=0, \quad x \in \RR^d
\end{equation}
with the property 
$$u_{c,A}(c)=A, \quad \|x\|^{d-\alpha}u_{c,A}(\|x\|) \to
K(d, \alpha, \beta)
\mbox { as }\|x\| \to \infty$$ where $K(d, \alpha,
\beta)$ is a positive constant. 
A natural conjecture now is that, like in the Laplacian case,
for $\beta = (d+\alpha)/(d-\alpha)-1$ the
$u_{c,A}$  constitute  {\em all} the bounded solutions of (\ref{ellbeta}),
and that for $\beta < (d+\alpha)/(d-\alpha)-1$ there are no bounded
non-zero solutions at all. 

As a first step to answer these questions, in Section \ref{last.section}
we make use of the  so called ``method of moving planes'', which is well known
in the Laplacian case, to show symmetry of positive solutions in a ball. Essential tools (like Hopf's boundary
lemma) can be carried over to the $\alpha$-Laplacian case. With this
method, we were able to show (see Theorem \ref{u.is.symmetric})
that the following equation
has only radially symmetric solutions:
\begin{eqnarray*}\label{ellB}
\Delta_{\alpha}u(x) + F(u(x))=0, \quad x \in B,\\
u(x) = 0, \quad x \in \RR^d \setminus B.
\end{eqnarray*}
Here $B$ is an open ball centered around $0$ and $F: \RR_{+} \to
\RR_{+}$ is non-decreasing.
We conjecture that an analogous statement is valid
for $F(z)=z^p$ with $p=(d+\alpha)/(d-\alpha)$ and $\RR^d$ instead
of $B$. In Section \ref{remarks} we describe several problems we think one
would have to overcome for carrying over the moving planes
method to the unbounded space setting in this case. 

\section{Two comparison criteria}
\setcounter{equation}{0}
\setcounter{lemma}{0}
\setcounter{proposition}{0}
\setcounter{corollary}{0}
\setcounter{theorem}{0}

In this section we assume that the function $G$ in Eq. \Ref{IVP}
satisfies the conditions
\begin{equation}
\label{G1} \lim_{z\to0+}\frac{G(z)}{z^{1+\beta}} = c\in(0,\infty)
\end{equation}
and
\begin{equation}
\label{G2} \int_{\theta}^{\infty}\frac{dz}{G(z)}<\infty
\end{equation}
for certain positive numbers $\beta$ and $\theta$.

\begin{lemma} \label{2.1} Let $G$ be a convex function satisfying \Ref{G1}, and
$\varepsilon>0$.
For any $M>0$ there exists $\varepsilon'>0$ such that
$$\frac{G((1+\varepsilon)z)}{(1+\varepsilon)z} >(1+ {\varepsilon'})
\frac{G(z)}{z} \mbox{ \ for $0<z\le M$. \ }
$$\label{lemma1}
\end{lemma}

{\noindent\bf Proof } By considering $G/c$ instead of $G$ we can assume that $c=1$. Given $\tilde{\varepsilon}>0$ there exists
$\delta>0$ such that
$$
(1-\tilde{\varepsilon})z^{\beta}<\frac{G(z)}{z}<(1+\tilde{\varepsilon})z^{\beta}
$$
for $z\in(0,\delta)$. Take $\tilde{\varepsilon}<
((1+\varepsilon)^{\beta} -1)/((1+\varepsilon)^{\beta}+1)$. Then
for $z<\delta/(1+\varepsilon):=x_0$,
$$
\frac{G((1+\varepsilon)z)}{(1+\varepsilon)z} >
(1-\tilde{\varepsilon})(1+\varepsilon)^{\beta}z^{\beta} >
\frac{(1-\tilde{\varepsilon})
(1+\varepsilon)^{\beta}}{(1+\tilde{\varepsilon})}
\frac{G(z)}{z}=(1+c')\frac{G(z)}{z}
$$
where $c'>0$. Since  $z\mapsto G(z)/z$ is continuous and strictly
increasing in $(0,\infty)$ it follows that
$\inf_{z\in[x_0,M]}\L(\frac{G((1+\varepsilon)z)}{(1+\varepsilon)z}/\frac{G(z)}{z}\R)
>1+c''$ with $c''>0$.
Taking $\varepsilon'=c' \wedge c''$ yields the assertion.
\hfill$\Box$
\medskip

Let us observe that if  $v$ is a
globally bounded solution of \Ref{IVP} we
necessarily have, for all $x$,
$$
P_tv(0,x)\to0\mbox{ \ as \ }t\to \infty,
$$
where $(P_t)$ is the semigroup with generator $L$. Indeed, from
the integral form of \Ref{IVP}
\begin{eqnarray*}
v(t,x) & = & P_tv(0,x) + \int_0^tP_{t-s}G(v(s,x))\,ds \\
& \ge & P_t v(0,x) + \int_0^tP_{t-s}G(P_sv(0,x))\,ds \\
& \ge & P_tv(0,x) + \int_0^tG(P_tv(0,x))\,ds \\
& \ge & tG(P_tv(0,x)),
\end{eqnarray*}
where we used in the first inequality that  $P_tv(0,\cdot)\le
v(t,\cdot)$, and Jensen's inequality after the second line. It follows from the global boundedness of $v$ that
\begin{equation}
\label{G3} \lim_{t\to\infty}P_tv(0,x)\le
\lim_{t\to\infty}G^{-1}(\mathrm{Const.} \, t^{-1})=0.
\end{equation}

\begin{proposition} 
\label{comp1}
Let $G$ be a convex, increasing function
satisfying \Ref{G1} and \Ref{G2}.
Assume the initial value $\varphi\ge 0$ leads to a globally bounded
solution of (\ref{IVP}). Then any initial value $\psi$ with  $0 \le \psi \le
(1-\varepsilon)\varphi$, $\varepsilon > 0$,
gives rise to a solution converging uniformly to zero.
\end{proposition}

{\noindent\bf Proof } Recall that the Feynman-Kac representation
of solutions of \Ref{IVP} is given by (see \cite{BLMW})
\begin{equation}
\label{FK4u}
u(t,x)=\int_{\RR^d}u(0,y)p_t(y,x)\E_y\L[\exp\L.\int_0^t\frac{G(u(s,X_s))}{u(s,X_s)}\,ds
\R|\,X_t=x \R]\,dy,
\end{equation}
where $(X_t)$ is the L\'{e}vy process with generator $L$, and
$p_t(x,y)$, $t>0$, $x,y\in\RR^d$, are its transition densities.

Suppose that $v$ is a globally bounded solution of
\Ref{IVP} and that $0\le  u(0,\cdot)\le(1-\varepsilon)v(0,\cdot)$
where $0<\varepsilon<1$. As \Ref{IVP} preserves ordering we have
$u(t,x)\le v(t,x)$ for all $t\ge0$ and $x\in \RR^d$, which
together with \Ref{FK4u} improves to
\begin{eqnarray*}
u(t,x) & \le & \int_{\RR^d}(1-\varepsilon)v(0,y)p_t(y,x)\E_y
\L[\exp\L.\int_0^t\frac{G(v(s,X_s))}{v(s,X_s)}\,ds
\R|\,X_t=x \R]\,dy\\
& = & (1-\varepsilon)v(t,x)
\end{eqnarray*}
uniformly in $t$ and $x$. Inserting this bound again into the
Feynman-Kac representation of $u$ yields
$$
u(t,x)\le\int_{\RR^d}u(0,y)p_t(y,x)\E_y\L[\exp\L.\int_0^t\frac{G((1-\varepsilon)
v(s,X_s))}{(1-\varepsilon)v(s,X_s)}\,ds \R|\,X_t=x \R]\,dy.
$$
Putting $z(t,x):= (1-\varepsilon)v(t,x)$ in the above inequality
and using Lemma \ref{lemma1} (with $\varepsilon$ in Lemma \ref{lemma1} 
substituted by
$\tilde{\varepsilon}:=\varepsilon/(1-\varepsilon)$) we get
\begin{eqnarray*}
u(t,x)&\le
&\int_{\RR^d}u(0,y)p_t(y,x)\E_y\L[\exp\L.\int_0^t\frac{1}{1+{\varepsilon}'}
\frac{G((1+\tilde{\varepsilon})z(s,X_s))}{(1+\tilde{\varepsilon})z(s,X_s)}\,ds
\R|\,X_t=x \R]\,dy \\
&=& \int_{\RR^d}u(0,y)p_t(y,x)\E_y\L[\L.e^{(1-\varepsilon'')A_t}
\R|\,X_t=x \R]\,dy,
\end{eqnarray*}
where  $\varepsilon'>0$ is given by Lemma \ref{lemma1}, \
$\varepsilon'':=\varepsilon'/(1+\varepsilon')$ and
\begin{equation}
\label{A_t} A_t:= \int_0^t\frac{G(v(s,X_s))}{v(s,X_s)}\,ds.
\end{equation}
Thus,
\begin{eqnarray*}
u(t,x) &\le & \int_{\RR^d}u(0,y)p_t(y,x)\E_y\L[ \L.
e^{(1-\varepsilon'')A_t}1(e^{A_t}<
\L(P_tv(0,x)\R)^{-1/2})\R|\,X_t=x\R]\,dy \\
& &  + \int_{\RR^d}u(0,y)p_t(y,x)\E_y\L[ \L. e^{
(1-\varepsilon'')A_t }1(e^{A_t}\ge
\L(P_tv(0,x)\R)^{-1/2})\R|\,X_t=x\R]\,dy.
\end{eqnarray*}
Since  $e^{A_t}\ge (P_tv(0,x))^{-1/2}$ implies
$e^{(1-\varepsilon'')A_t} \le e^{A_t}
\L(P_tv(0,x)\R)^{\varepsilon''/2}$, we obtain
\begin{eqnarray*}
u(t,x)& \le & \int_{\RR^d}u(0,y)p_t(y,x)\L(P_tv(0,x)\R)^{-1/2} \\
& & +
\L(P_tv(0,x)\R)^{\varepsilon''/2}\int_{\RR^d}u(0,y)p_t(y,x)\E_y\L[\L.e^{A_t}\R|\,X_t=x
\R]\,dy \\
& \le &  (1-\varepsilon)\L((P_tv(0,x))^{1/2} +
v(t,x)(P_tv(0,x))^{\varepsilon''/2} \R)
\end{eqnarray*}
which tends to $0$ uniformly as $t\to\infty$ due to \Ref{G3}.
\hfill $\Box$

\begin{proposition} \label{comp2}  Let $G$ be a convex, increasing function
satisfying \Ref{G1} and \Ref{G2}.
Assume the initial value $\varphi=v(0,\cdot)\ge 0$ leads to a 
globally bounded 
solution of (\ref{IVP}) which for some open ball $B \subset \RR^d$  and some $\kappa>0$  obeys
\begin{equation}
\label{odd.condition} \inf_{x\in B}v(t,x)\ge \kappa \mbox{ \ for
all sufficiently large $t>0$.}
\end{equation}
Then for any $\varepsilon > 0$, the initial condition
$(1+\varepsilon)\varphi$
leads to blow-up in finite time.
\end{proposition}

{\noindent \bf Proof\ } 
    By the Feynman-Kac formula,
$$
v(t,x)=\int_{\RR^d}v(0,y)p_t(y,x)\E_y\L[\L.
e^{A_t}\R|\,X_t=x\R]\,dy,
$$
where $A_t$ is given by \Ref{A_t}. If $K>0$ then
$$
\int_{\RR^d}v(0,y)p_t(y,x)\E_y\L[\L.e^{A_t}; A_t\le
K\R|\,X_t=x\R]\,dy \le e^K\E_xv(0,X_t) \to 0
$$
as $t\to\infty$ uniformly in $x$ due to \Ref{G3}. Therefore, for
all $K>0$ 
there exists $T_0=T_0(K,\gamma)>0$ such
that for $t> T_0$
\begin{equation}
\label{upper.bound.t}
\frac{v(t,x)}{\int_{\RR^d}v(0,y)p_t(y,x)\E_y\L[\L. e^{A_t}; A_t\ge
K\R|\,X_t=x\R]\,dy} \le 2. 
\end{equation}
Without loss of generality we can assume that $\inf_{x\in
B_1(0)}v(t,x)\ge\kappa$ for all $t$ large enough, where $B_r(x)$
denotes the ball in $\RR^d$ of radius $r$ centered at $x$. Arguing
as above we check via the Feynman-Kac representation that for all
$t\ge0$ and $x\in \RR^d$, $u(t,x)\ge (1+\varepsilon)v(t,x)$. Plugging this again into the
Feynman-Kac representation for $u$ yields
\begin{eqnarray*}
u(t,x) & \ge &
(1+\varepsilon)\int_{\RR^d}v(0,y)p_t(y,x)\E_y\L[\L.\exp\int_0^t\frac{G((1+\varepsilon)
v(s,X_s))}
{(1+\varepsilon)v(s,X_s)}\,ds\R|\,X_t=x\R]\,dy \\
& \ge &
(1+\varepsilon)\int_{\RR^d}v(0,y)p_t(y,x)\E_y\L[\L.e^{(1+\varepsilon')A_t}
\R|\,X_t=x\R]\,dy
\end{eqnarray*}
for some $\varepsilon'>0$ by Lemma \ref{2.1}. Using this and
\Ref{upper.bound.t} we obtain for given $K>0$ 
and $t$
big enough that
\begin{eqnarray*}
u(t,x) &\ge &
(1+\varepsilon)e^{K\varepsilon'}\int_{\RR^d}v(0,y)p_t(y,x)\E_y\L[\L.
e^{A_t}; A_t\ge K\R|\,X_t=x\R]\,dy \\
& \ge & (1+\varepsilon)e^{K\varepsilon'}\frac{v(t,x)}{2}.
\end{eqnarray*}
Hence, for any $K > 0$ we find
$\inf_{x\in B_1(0)} u(t,x) \ge \kappa 
(1+\varepsilon)e^{K\varepsilon'}/2$ 
for all
sufficiently large $t$. As is well known (see e.g. \cite{KST}),
this inequality together with \Ref{G2} are sufficient for
finite-time blowup of $u$. \hfill$\Box$

\begin{corollary} Let $G$ be a convex, increasing function
satisfying \Ref{G1} and \Ref{G2}, and $\varphi \ge 0$  a non-trivial positive bounded solution of
$$
L\varphi(x) + G(\varphi(x))=0.
$$

a) For each $\varepsilon>0$ the solution $u$ of \Ref{IVP} with
initial value $u(0,x) = (1+\varepsilon)\varphi(x)$, $x\in\RR^d$, blows
up in finite time.

b) For each  $\varepsilon\in(0,1)$ the solution $u$ of \Ref{IVP} with
initial value $u(0,x) = (1-\varepsilon)\varphi(x)$, $x\in\RR^d$ converges 
uniformly to 0 as $t\to\infty$.

\end{corollary}

{\noindent\bf Proof } This is immediate from propositions \ref{comp1} 
and \ref{comp2}.
\hfill$\Box$

\section{A class of radially symmetric stationary solutions}
\setcounter{equation}{0}
\setcounter{lemma}{0}
\setcounter{proposition}{0}
\setcounter{corollary}{0}
\setcounter{theorem}{0}

We now set out to specify a family of positive stationary solutions of
\Ref{IVP} in the particular case of
$G(z)=z^p$, where $p=(d+\alpha)/(d-\alpha)$. Before doing this, we
still consider the case of a general $p$, and note that
the ``elliptic'' equation
\begin{equation}\label{elliptic}
\Delta_{\alpha}u(x) + u^p(x)=0, \ \ x\in\RR^d
\end{equation}
can be rewritten in integral form as
\begin{equation}\label{integral}
     u(x)= \int_0^{\infty}\E_x\L[u^p(X_t)\R]\,dt, \ \ x\in\RR^d,
\end{equation}
where  $(X_t)$ denotes the (symmetric) $\alpha$-stable process in
$\RR^d$. Hence, for $d\le\alpha$, due to recurrence of $(X_t)$, the only
non-negative solutions of (\ref{elliptic}) are $u \equiv 0$ and
$u\equiv\infty$.
Therefore, we henceforth assume that $d>\alpha$,
in which case (\ref{integral}) rewrites as 
(see \cite{BG}, p. 264)
\begin{equation}\label{Green}
u(x)= \int_{\RR^d}\frac{\AA u^p(y)}{\|y-x\|^{d-\alpha}}\,dy, \ \ x\in\RR^d, 
\end{equation}
where $\AA:=\Gamma(\frac{1}{2}(d-\alpha))/[\Gamma(\frac{1}{2}\alpha)2^{\alpha}\pi^{d/2}]$.

\begin{proposition}\label{steadystates}
If $p=(d+\alpha)/(d-\alpha)$, then for any $A\in(0,\infty) $ and $c\in \RR^d$ the function
$$u_{c,A}(x)=\frac{A}{\L[1+\L(A^{2/(d-\alpha)}2^{-1}\L(
\Gamma\L(\frac{d+\alpha}{2}\R) \L/
\Gamma\L(\frac{d-\alpha}{2}\R)\R.
\R)^{-1/\alpha}\|x-c\| \R)^2\R]^{(d-\alpha)/2}}
$$
solves \Ref{elliptic}.%
\end{proposition}

{\noindent \bf Proof } Without loss of generality we assume that $c$ is the origin.
 Due to \Ref{Green} it suffices to show that
\begin{equation}\label{claim}
u_{0,A}(x) =\int_{\RR^d}\frac{\AA u_{0,A}^p(y)}{\|y-x\|^{d-\alpha}}\,dy, \ \ x\in\RR^d.
\end{equation} 
Let us write $a:= A^{-2/(d-\alpha)}2\L(
\Gamma\L(\frac{d+\alpha}{2}\R) \L/
\Gamma\L(\frac{d-\alpha}{2}\R)\R.
\R)^{1/\alpha} $. We first note that 
$$
u_{0,A}(x) = \frac{A}{(1+4\pi^2\|\frac{x}{2\pi a}\|^2)^{(d-\alpha)/2}}
= A\widehat{B_{d-\alpha}}\L(\frac{x}{2\pi a}\R) $$
(see \cite{Folland}, p. 155), where for any $f\in L^1(\RR^n)$, $\widehat{f}(x):=\int_{\RR^d}e^{-2\pi iy\cdot x}f(y)\,dy$ is the Fourier transform of $f$, and for any complex $w$ with $\mathrm{Re}(w)>0$
$$
B_w(x):= \frac{1}{\Gamma(\frac{w}{2})(4\pi)^{d/2}}\int_0^{\infty}r^{(w-d)/2-1}e^{-r-\|x\|^2/4r}
\,dr.
$$ 
Hence
$$
\widehat{u_{0,A}}(x)  =  A\L[\widehat{B_{d-\alpha}}\L(\frac{\cdot}{2\pi a}\R)  \R]\!\widehat{\phantom{|}}(x) =  A(2\pi a)^d\widehat{\widehat{B_{d-\alpha}}}(2\pi a x) 
 = A(2\pi a)^d B_{d-\alpha}(-2\pi a x).
$$
Notice that 
\begin{equation}
\label{B.M}
B_w(x)= \frac{2^{(d-w)/2+1}}{\Gamma\L(\frac{w}{2}\R)(4\pi)^{d/2}}     \|x\|^{(w-d)/2}K_{(d-w)/2}(\|x\|), \ \ x\in \RR^d,
\end{equation}
where for any complex $\nu$, 
$K_{\nu}$ is the Macdonald's function 
(\cite{Watson}, \S 6$\cdot$22), also known as modified Bessel 
function of the second kind, which is given by
$$K_{\nu}(z)=\frac{1}{2}\L(\frac{1}{2}z\R)^{\nu}\int_0^{\infty}r^{-\nu-1}e^{
-r-z^2/4r}dr,\ \ \mathrm{Re}(z^2)>0.
$$ 
 It follows that
\begin{equation}
\label{LHS.transform}
\widehat{u_{0,A}}(x) = A a^{d-\alpha/2}\pi^{(d-\alpha)/2}\frac{2}{\Gamma(\frac{d-\alpha}{2})}\,\|x\|^{-\alpha/2}
K_{\alpha/2}(\|
2\pi a x\|).
\end{equation}  
To compute the Fourier transform of the other side of \Ref{claim} we use the convolution theorem, \Ref{B.M} and 
$$\L[\AA\|\cdot\|^{-(d-\alpha)}\R]\!\widehat{\phantom{|}}(x) = \L(2\pi\|x\|\R)^{-\alpha},\ \ \ x\in\RR^d,\ \ \
0<\mathrm{Re}(\alpha)<d,
$$
(e.g. \cite{Folland}, p. 154)
to obtain 
\begin{eqnarray*}
\L[\int_{\RR^d}\!\!\frac{\AA u_{0,A}^p(y)}{\|y-\cdot\|^{d-\alpha}}\,dy\R]\!\widehat{\phantom{|}}\,(x)
& = & (2\pi)^{-\alpha}\|x\|^{-\alpha}A^p\L[\frac{1}{\L(1+4\pi^2\|\frac{\cdot}{2\pi a}\|^2\R)^{(d+\alpha)/2}}
\R]\!\widehat{\phantom{|}}\,(x) \\
& = & (2\pi)^{-\alpha}\|x\|^{-\alpha}A^p
\L[\widehat{B_{d+\alpha}}\L(\frac{\cdot}{2\pi a}\R)  \R]\!\widehat{\phantom{|}}(x)\\
& = & (2\pi)^{-\alpha}\|x\|^{-\alpha}A^p(2\pi a)^dB_{d+\alpha}(-2\pi ax)\\
& = & 2^{-\alpha}\pi^{(d-\alpha)/2 }\|x\|^{-\alpha/2}A^pa^{d+\alpha/2}\frac{2}{\Gamma\L(\frac{d+\alpha}{2}\R)}
K_{-\alpha/2}(\|2\pi ax\|)
\end{eqnarray*}
Since $K_{\nu}=K_{-\nu}$ (\cite{AS}, Formula 9.6.6), by comparing
the RHS of the last equality with that of \Ref{LHS.transform} we see that they are equal for the value of $a$ stated at the beginning of the proof. The result follows from uniqueness of Fourier transforms.\hfill$\Box$
\medskip

{\noindent\bf Remarks\ } 1.
Recall \cite{BB,GT} that for $0<\alpha\le
2$ and $d>\alpha$ the Kelvin transform of $u$ is defined by
\begin{equation}\label{Kelvin1}
v(x):=\frac{1}{\|x\|^{d-\alpha}}u\L(\frac{x}{\|x\|^2}\R),\ \ x\in
\RR^d, \ \ x\neq0.
\end{equation}
A simple calculation shows that for any fixed $c\in\RR^d$ the family of solutions $\{u_{c,A}\}_{A\ge0}$ rendered by Proposition \ref{steadystates} is invariant under the Kelvin transform  with center  $c$
$$v_c(x) := \frac{1}{\|x-c\|^{d-\alpha}}u_{c,A}\L(c+\frac{x-c}{\|x-c\|^2}\R).
$$

2. Moreover, if $u$ is {\em any\/} given regular positive solution of \Ref{elliptic} (where $0<\alpha\le 2$, $p>1$ and $d>\alpha$) and $x\neq0$, then its Kelvin transform $v(x)$
satisfies \begin{equation}\label{Kelvin2}
\Delta_{\alpha} v(x)+\frac{v^p(x)}{\|x\|^{(d+\alpha)
-p(d-\alpha)}}=0.
\end{equation} 
Indeed, let $G_{\alpha}$ denote the Green's
operator corresponding to $\Delta_{\alpha}$, and  $x\neq0$. Then,
\begin{eqnarray*}
-G_{\alpha}\L(-\frac{v^p(x)}{\|x\|^{(d+\alpha) -p(d-\alpha)}} \R)
& = &
\int\frac{\AA}{\|x-y\|^{d-\alpha}}\cdot\frac{1}{\|y\|^{(d+\alpha)-p(d-\alpha)}}\\ & & \hspace{7.5em} \cdot
\frac{1}{\|y\|^{(d-\alpha)p}}\cdot u^p\L(\frac{y}{\|y\|^2}\R)\,dy\\
& = & \AA\int \frac{1} {\L\|x-\frac{z}{\|z\|^2}\R\|^{d-\alpha}}
\cdot \frac{1}{\|z\|^{d-\alpha}}\cdot u^p(z)\,dz \\
&  =& \AA\int\frac{1}{\L\| \:x\cdot\|z\| - \frac{z}{\|z\|}
\R\|^{d-\alpha}}\cdot u^p(z)\,dz\\
& = &\AA\int_{\RR^d}\frac{u^p(z)\,dz}{\L\|
\:\frac{x}{\|x\|}-\|x\|\cdot z\:
\R\|^{d-\alpha}}\\
&=&\frac{1}{\|x\|^{d-\alpha}}u\L(\frac{x}{\|x\|^2}\R) = v(x),
\end{eqnarray*}
where we used the elementary identity $\|\: x\cdot\|z\| -\frac{z}{\|z\|}\| =
\|\:\frac{x}{\|x\|}-\|x\|\cdot z\|$ in the fourth equality.

3. Proceeding as in the proof of Proposition \ref{steadystates} one can verify
that a singular explicit solution
to \Ref{elliptic} is given by 
$$u_{\mathrm{sing}}(x)= {\textstyle \L[2^{\alpha} \L(\Gamma\L(\frac{d+\alpha}{4}\R)\L/
\Gamma\L(\frac{d-\alpha}{4}\R.\R)  \R)^2 \R]^{1/(p-1)}}\cdot\:\frac{1}{\|x\|^{(d-\alpha)/2}},\ \ \ x\neq0,$$  and that $u_{\mathrm{sing}} $ is a fixed point of the
Kelvin transform \Ref{Kelvin1}.

\section{Rotational symmetry of solutions in a ball}\label{last.section}
\setcounter{equation}{0}
\setcounter{lemma}{0}
\setcounter{proposition}{0}
\setcounter{corollary}{0}
\setcounter{theorem}{0}

Let $F : \RR_+ \to \RR_+$ be non-decreasing, not identically constant, and
$\alpha \in (0,2)$.
Let $u : \RR^d \to \RR_+$ be a non-negative bounded solution to
\begin{equation}
\label{steadystate}
\Delta_\alpha u + F(u) = 0, \; x \in B_1(0), \quad
        u \equiv 0 \:\; \mbox{in $B_1(0)^c$},
\end{equation}
i.e., for all $x$ in the unit ball we have
\begin{equation}
\label{steadystate.intform}
u(x) = 
\int_0^\infty \E_x\left[ F(X_t); \,{\textstyle \sup_{u \le t} |X_u| < 1}
\right] dt
= \int_{B_1(0)} G_{\alpha}(x,y) F(u(y)) \,dy,
\end{equation}
where $X$ is the symmetric $\alpha$-stable process and 
$G_{\alpha}(x,y)$ is the corresponding Green's function for the unit ball. 

That \Ref{steadystate} has non-trivial positive solutions follows from \cite{BN}
at least in the case $\alpha=2$. Without loss of generality we assume that $u \not\equiv 0$.
Our aim here is to show that
\begin{theorem}\label{u.is.symmetric} $u$ is rotationally symmetric about the origin.
\end{theorem}
Our approach is based on the {\em method of moving planes,\/}  
a device that goes back to Alexandrov \cite{Alex} and
has by now a venerable history in the study of symmetries of
solutions of pde's, see also \cite{Serrin} and \cite{GNN}. 
The idea is as follows:
Choose any direction in $\RR^d$, wlog the $x_1$-direction,
and show that $u$ is mirror symmetric with respect to the
hyperplane through the origin with this given direction as a
normal vector.
In order to achieve this let us define for $\lambda \in (-1,1)$
\[
T_\lambda := \{ x \in \RR^d : x_1 = \lambda \},
\quad
\Sigma_\lambda := \{ x \in \RR^d : x_1 < \lambda \}
\]
and for $x \in \RR^d$ let $x^\lambda := (2\lambda-x_1,x_2,\ldots,x_d)$
be the image under reflection along $T_\lambda$;
see Figure 1.
\begin{figure}
\begin{center}
\psfrag{x}{\colorbox{white}{$x$}}
\psfrag{0}{$0$}
\psfrag{1}{$1$}
\psfrag{-1}{\colorbox{white}{$-1$}}
\psfrag{x1}{$x_1$}
\psfrag{xlambda}{$x^\lambda$}
\psfrag{Slambda}{\colorbox{white}{$\Sigma_\lambda$}}
\psfrag{Tlambda}{$T_\lambda$}
\psfrag{lambda}{$\lambda$}
\includegraphics{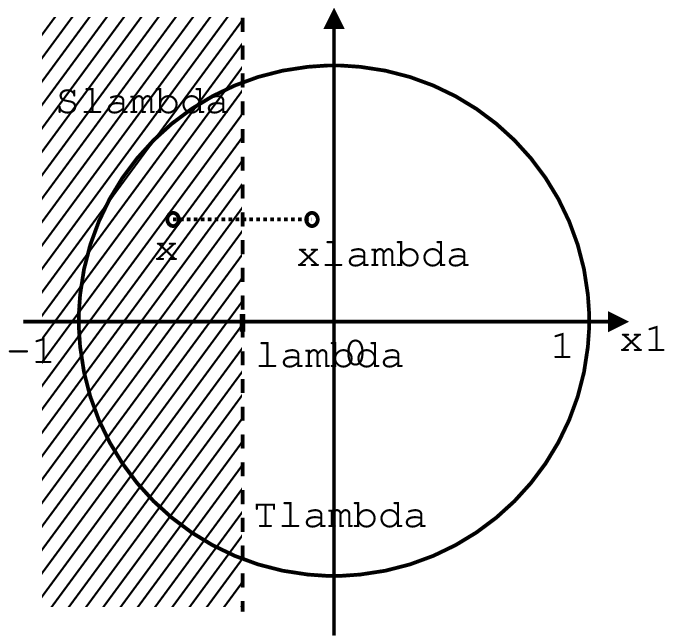}
\end{center}
\centerline{Figure 1}
\end{figure}

Define the set $\Lambda$ by
\begin{equation}
\Lambda := \left\{ \lambda \in (-1,0) : u(x^\lambda) \ge u(x) \;
\forall \, x \in \Sigma_\lambda,
\frac{\partial}{\partial x_1} u(x)>0 \;
\forall \, x \in T_\lambda \cap B_1(0) \right\}.
\end{equation}
Observe that by the minimum
principle we have $u(x^\lambda) > u(x)$
for $x \in \Sigma_\lambda \cap B_1(0)$ and $\lambda \in \Lambda$,
i.e.,  $\lambda \in \Lambda$ means that reflection along $T_\lambda$
(strictly) increases the value of $u$.
In Section \ref{proof.that.zero.is.okay} we prove that
\begin{equation}
\label{zero.is.okay}
\sup \Lambda = 0
\end{equation}
so that by continuity
$u(-x_1,x_2\ldots,x_d) \ge u(x_1,x_2\ldots,x_d)$ whenever $x_1\le 0$.
By considering $\lambda>0$ and working in the opposite direction
we can then conclude the reversed inequality and hence obtain the
desired symmetry.

\subsection{Some preparatory lemmas}
\begin{lemma}
A bounded solution $u$ of (\ref{steadystate.intform})
satisfies $u \in C(\RR^d) \cap C^\infty(B_1(0))$.
\end{lemma}
{\noindent\bf Proof \ } Using the explicit form of the
Green's kernel of the ball (see e.g. \cite{BB}, formula (2.3))
one easily checks that interchange of integration and
differentiation is justified.

\begin{lemma}[Minimum principle] Let $D\subset\RR^d$ be a bounded domain. Suppose 
$u : \RR^d \to \RR_+$ is continuous, $\Delta_\alpha u \le 0$ on $D$, and satisfies
$u \equiv 0$ on $D^c$.
Then either $u \equiv 0$ or $u > 0$ on $D$.\label{min.pri}
\end{lemma}
{\noindent\bf Proof\ }
Let $D_\varepsilon :=\{ x \in D : u(x)> \varepsilon \}$.
By continuity $D_\varepsilon$ is open.
Assume that $D_\varepsilon \neq \emptyset$ for some
$\varepsilon >0$.
Let $(X_t)$ be the $\alpha$-stable process, and
$\tau :=\inf \{s : X_s \not\in D\}$ the hitting time of $D^c$.
Then $M_t := u(X_{t \wedge \tau}) - u(X_0) -
\int_0^{t \wedge \tau} \Delta_\alpha u(X_s) \, ds$
is a $\Pr_x$-martingale for each $x \in D$.
Let furthermore $\tau' := \inf \{s : X_s \in D_\varepsilon \}$.
For each $x \in D$ we have $\E_x M_{\tau'} = 0$, or
\[
u(x) = \E_x\left[ u(X_{\tau' \wedge \tau})\right] +
\E_x\left[ \int_0^{\tau' \wedge \tau} (-\Delta_\alpha u)(X_s) \,ds \right]
\ge \E_x\left[ u(X_{\tau' \wedge \tau})\right] \ge \varepsilon \Pr_x(\tau' < \tau) > 0
\]
because $(X_t)$ hits any open subset of $D$ with positive probability
before exiting from $D$. \hfill$\Box$

We will have occasion to consider the behavior of $u$ at the boundary
of the ball. In this respect, the following lemma is helpful:
\begin{lemma}[Hopf's $\alpha$-stable boundary lemma]\label{Hopf.lemma}
Let $D \subset \RR^d$ be open, $u : \RR^d \to \RR_+$ continuous
with $u \equiv 0$ on $D^c$, $\Delta_\alpha u \le 0$ on $D$, 
$u$ not identically zero. 
Let $x_0 \in \partial D$ satisfy an {\em interior sphere condition},
i.e. there exists a ball $B_\delta(x_1) \subset D$ with
$\overline{B_\delta(x_1)} \cap D^c = \{x_0\}$, and let
$\nu$ be an outward pointing unit vector at $x_0$.
Then
\[
\frac{\partial}{\partial \nu}u(x_0) < 0
\]
(in fact,
$\lim_{\varepsilon \searrow 0} (u(x_0)-u(x_0-\varepsilon \nu))/\varepsilon
= - \infty$).
\end{lemma}
{\noindent\bf Proof \ }
Because of the interior sphere condition at $x_0$ we can
find a ball $B_\delta(x_1) \subset D$ such that
$\overline{B_\delta(x_1)} \cap D^c = \{x_0\}$ and
also $\widetilde{B} \subset D$, where $\widetilde{B}$ is the
``left half'' of a spherical shell around $x_1$ with
interior radius $\delta$ and exterior radius $\delta'>\delta$; see Figure 2.
Observe that $u>0$ in $D$ by Lemma \ref{min.pri}, in particular $\inf_{\widetilde{B}} u >0$
because $\widetilde{B}$ is compact and $u$ continuous.
\begin{figure}
\begin{center}
\psfrag{btilde}{$\widetilde{B}$}
\psfrag{x0}{$x_0$}
\psfrag{x1}{$x_1$}
\psfrag{x}{$x$}
\psfrag{nu}{\colorbox{white}{$\nu$}}
\psfrag{d}{$D$}
\psfrag{dcomp}{\colorbox{white}{$D^c$}}
\psfrag{GAMMA}{$\gamma$}
\includegraphics{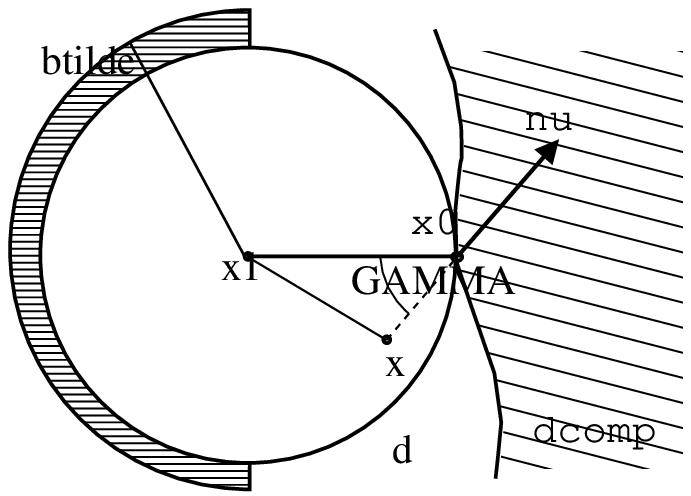}
\end{center}
\centerline{Figure 2}
\end{figure}
Let $\tau:=\inf\{ t : X_t \not\in B_\delta(x_1) \}$. Then
\[
u(x) = \E_x\left[ u(X_\tau) - \int_0^\tau (\Delta_\alpha u)(X_t) \,dt \right]
        \ge \E_x \, u(X_\tau)
\]
for $x\in B_{\delta}(x_1)$. Take $x=x_0-\varepsilon \nu$ with $\varepsilon>0$ small enough, and denote by $\gamma$ the angle between
$\overline{x_1\,x_0}$ and $\overline{x\,x_0}$.
Then $\gamma \in (-\pi/2, \pi/2)$ because $\nu$ is an
outward pointing vector,  hence $\cos \gamma > 0$. The cosine theorem gives
$|x-x_1|^2 = \delta^2 + \varepsilon^2 - 2\delta \varepsilon \cos \gamma$.
Using the explicit form of the Poisson kernel for the
complement of a ball (see e.g. \cite{BB}, formula (2.2) 
or \cite{BGR}, Thm. A and rescale) we can estimate
\begin{eqnarray*}
u(x) & \ge & \E_x \, u(X_\tau) \ge \int_{\widetilde{B}} P(x,y) u(y) \,dy \\
& \ge & \left( \inf_{\widetilde{B}} u \right)
        C_{\alpha, d} \int_{\widetilde{B}} \left(
                \frac{\delta^2-|x-x_1|^2}{|y-x_1|^2-\delta^2}
                                                \right)^{\alpha/2}
                |x-y|^{-d} \,dy \\
& \ge & C \left(\delta^2-|x-x_1|^2\right)^{\alpha/2}
\ge C' \varepsilon^{\alpha/2}.
\end{eqnarray*}
Thus we see that
$\limsup_{\varepsilon \searrow 0} (u(x_0)-u(x_0-\varepsilon \nu))/\varepsilon
= - \infty$.\hfill$\Box$

\begin{lemma}
\label{lemma.symmetric.boundary}
Let $w: \RR^d \to \RR$ be continuous, bounded with
$w \ge 0$ on $\Sigma_0$, $w(x^0) = -w(x)$.
Let $x_* \in T_0$ be such that there exists a $\delta >0$
with $\Delta_\alpha w \le 0$ on $B_\delta(x_*) \cap \Sigma_0$.
Then either
\[ w \equiv 0,  \mbox{\ \ or \ \ }
 w > 0 \; \mbox{on $B_\delta(x_*) \cap \Sigma_0$ and} \;
\frac{\partial}{\partial x_1}w(x_*) < 0. \]
\end{lemma}

{\noindent\bf Proof\ }
Assume $w \not\equiv 0$.
Let $x \in B_\delta(x_*) \cap \Sigma_0$.
If $w(x) = 0$ we would have (see e.g. \cite{BB})
\[
\Delta_\alpha w(x) = c_{\alpha,d}
        {\rm PV} \int_{\RR^d} \frac{w(x+y)-0}{|y|^{d+\alpha}} dy > 0
\]
by the non-triviality and symmetry of $w$, in contradiction to the
assumption.

To show that the derivative is non-zero choose $\delta'>0$
such that $\sup_{\Sigma_0 \setminus B_{\delta'}(x_*)} w >0$.
Let $(X_t)$ be the $\alpha$-stable process,
$\tau := \inf \{ t : X_t \not\in B_{\delta'}(x_*) \}$. Define
$v(x):= \E_x w(X_\tau)$,
$\widetilde{v}(x) := \E_x \int_0^\tau (-\Delta_\alpha w)(X_s) \, ds$.
Observe that
\[
 \Delta_\alpha v = 0 \: \mbox{ in $B_{\delta'}(x_*)$ and } \:
 v = w \: \mbox{ in $\RR^d \setminus B_{\delta'}(x_*)$},
\]
\[
 \Delta_\alpha \widetilde{v} = \Delta_\alpha w \:
        \mbox{ in $B_{\delta'}(x_*)$ and } \:
  \widetilde{v} = 0 \: \mbox{ in $\RR^d \setminus B_{\delta'}(x_*)$}.
\]
Uniqueness of the Dirichlet problem for $\Delta_\alpha$ in
$B_{\delta'}(x_*)$ thus gives $w = v + \widetilde{v}$.
We have $v(x)=\int_{|y-x_*|>\delta'} P(x,y) w(y) \,dy$ for
$x \in B_{\delta'}(x_*)$, where the Poisson kernel is given by
(see e.g. \cite{BB}, formula (2.2))
\[
P(x,y) = C_{\alpha,d} \left[ \frac{\delta'^2 - |x-x_*|^2}{%
                                |y-x_*|^2-\delta'^2} \right]^{\alpha/2}
                |x-y|^{-d}, \quad x \in B_{\delta'}(x_*),
                        y \not\in B_{\delta'}(x_*).
\]
One checks that $\frac{\partial}{\partial x_1}P(x_*,y) < 0$ for
$y \in \Sigma_0$ and $\frac{\partial}{\partial x_1}P(x_*,y) > 0$
for $y \in - \Sigma_0$. The interchange of integration and
differentiation is justified because $w$ is bounded, so we can
compute
\[
\frac{\partial}{\partial x_1} v(x_*) =
        \int_{|y-x_*|>\delta'} \frac{\partial}{\partial x_1}P(x_*,y) w(y) \,dy
        < 0.
\]
Furthermore, for $x \in B_{\delta'}(x_*)$
\[
\widetilde{v}(x) = \int_{B_{\delta'}(x_*)} G(x,y) (-\Delta_\alpha w)(y) \,dy,
\]
where the Green kernel for $B_{\delta'}(x_*)$ is given by
(see e.g. \cite{BB}, formula (2.3) or \cite{BGR} and consider
the obvious scaling properties of $\alpha$-stable processes)
\begin{equation}
\label{green.function.ball}
G(x,y) = c_{\alpha,d} (|x-y|)^{\alpha-d}
\int_0^{w_{\delta'}(x,y)} \frac{r^{\alpha/2-1}}{(r+1)^{d/2}} dr,
\end{equation}
where
$w_{\delta'}(x,y)=(\delta'^2-|x-x_*|^2)(\delta'^2-|y-x_*|^2)/|x-y|^2$.
Inspection shows that for $x, y \in B_{\delta'}(x_*) \cap \Sigma_0$
we have $G(x,y) \ge G(x,y^0)$, hence
$\widetilde{v} \ge 0$ in $B_{\delta'}(x_*) \cap \Sigma_0$
and by symmetry $\widetilde{v} \le 0$
in $B_{\delta'}(x_*) \cap (-\Sigma_0)$.
We conclude that $(\partial/\partial x_1)\widetilde{v}(x_*) \le 0$ and thus
\[
\frac{\partial}{\partial x_1} w(x_*) \le
        \frac{\partial}{\partial x_1} v(x_*) < 0.
\]\hfill$\Box$
\subsection{Proof of (\ref{zero.is.okay})}
\label{proof.that.zero.is.okay}

We proceed in three steps and show that
\begin{enumerate}
\item $\Lambda \supset (-1,-1+\varepsilon)$ for some 
small $\varepsilon > 0$,
\item $\Lambda$ is open. In particular for $\lambda \in \Lambda$
there exists $\varepsilon > 0$ such that
$[\lambda, \lambda + \varepsilon) \subset \Lambda$.
\item From 1. and 2. we conclude that
$\Lambda = (-1,\lambda_{\rm max})$. We finally show
that $\lambda_{\rm max}=0$.
\end{enumerate}
%

{\noindent
\em Step 1.\/ }
Obviously $(-1,0,\ldots,0)$ is an outward pointing direction
for each $x \in \partial B_1(0) \cap \Sigma_{-1/2}$.
By the boundary lemma and the fact that $(\partial/\partial x_1)u$
is continuous in $B_1(0)$  there is an open neighborhood
$D$ of $\partial B_1(0) \cap \Sigma_{-1/2}$ such that
$(\partial/\partial x_1)u > 0$ on $D \cap B_1(0)$.
Choose $\varepsilon > 0$ so small that
$\Sigma_{-1+\varepsilon} \cap B_1(0) \subset D$.
Then $-1+\varepsilon/2 \in \Lambda$.
\\[2ex]
{\noindent\em Step 2.\/} 
We argue by contradiction. Assume there was $\lambda_* \in \Lambda$
and also a sequence $(\lambda_n) \subset (-1,0) \setminus \Lambda$
with $\lambda_n \searrow \lambda_*$.
>From the definition of $\Lambda$, possibly passing to a
suitable subsequence (which we again would denote by $(\lambda_n)$)
we can always arrive at one of the following possibilities:
\begin{itemize}
\item[a)] There exists a sequence $(x_n) \subset B_1(0)$,
$x_n \in \Sigma_{\lambda_n}$, with
$x_n \to x_* \in \overline{B_1(0)}$ and $u(x_n) \ge u(x_n^{\lambda_n})$
for all $n$, or
\item[b)] There exists a sequence $(x_n) \subset B_1(0)$,
$x_n \in T_{\lambda_n}$,
with
$x_n \to x_* \in \overline{B_1(0)}$ and
$(\partial/\partial x_1)u(x_n) \le 0$ for all $n$.
\end{itemize}
Assume a) was true. We cannot have $x_* \in \Sigma_{\lambda_*}$
because $u$ is continuous and $u(x^{\lambda_*}) > u(x)$
for $x \in \Sigma_{\lambda_*}$ by the above remark.
Hence $x_* \in T_{\lambda_*} \cap \overline{B_1(0)}$.
But then we have $(\partial/\partial x_1)u(x_*) =
\lim_{n\to\infty} (u(x_n^{\lambda_n})-u(x_n))/(2{\rm d}(x_n,T_{\lambda_n})) \le 0$.
By Hopf's boundary lemma, this forces $x_*$ to be away from
$\partial B_1(0)$, but then we obtain a contradiction to
$\lambda_* \in \Lambda$.

If b) was true we would again find a point
$x_* \in T_{\lambda_*} \cap \overline{B_1(0)}$
with $(\partial/\partial x_1)u(x_*) \le 0$ and arrive at a
contradiction.
\\[2ex]
{\noindent\em Step 3.\/} 
>From the preceding steps we know that $\Lambda = (-1, \lambda_{\rm max})$.
If $u(x^{\lambda_{\rm max}}) \equiv u(x)$ for
$x \in \Sigma_{\lambda_{\rm max}}$ then we have found a symmetry
center. As $u$ is continuous, $u=0$ on $\partial B_1(0)$ and strictly
positive inside $B_1(0)$ this can only be true for
$\lambda_{\rm max}=0$. Indeed, if  $\lambda_{\rm max} < 0$ then by continuity we would have
$u(x) \le u(x^{\lambda_{\rm max}})$ for $x \in \Sigma_{\lambda_{\rm max}}$,
but with $u(x) \not\equiv u(x^{\lambda_{\rm max}})$.
Define $w(x) := u(x^{\lambda_{\rm max}}) - u(x)$. Observe that
$w$ is continuous and bounded, non-negative in  $\Sigma_{\lambda_{\rm max}}$ and
$w(x^{\lambda_{\rm max}}) = - w(x)$.
For $x \in \Sigma_{\lambda_{\rm max}} \cap B_1(0)$ we have
\[
\Delta_\alpha w(x) = (\Delta_\alpha) u(x^{\lambda_{\rm max}})
- (\Delta_\alpha) u(x) = -\L[F (u (x^{\lambda_{\rm max}} ) ) - F(u(x))\R]
\le 0,
\]
and we infer from Lemma \ref{lemma.symmetric.boundary} that
$(\partial/\partial x_1)w(x) < 0$ for all
$x \in T_{\lambda_{\rm max}} \cap B_1(0)$. 
In conclusion,  $\lambda_{\rm max} < 0$  implies 
$\lambda_{\rm max} \in \Lambda$ which by Step 2
forces $\sup \Lambda > \lambda_{\rm max}$.
This is a contradiction.

\subsection{Remarks and open questions}\label{remarks}

The approach developed in the last sections invites to try to use 
the moving planes method also for solutions $u$ of (3.1), 
that is to show radial
symmetry about some point for the corresponding problem
on all of $\RR^d$. In fact, this has been done in the case
of the classical Laplacian starting with the work of Serrin \cite{Serrin}; 
see also \cite{Bianchi,CL,GNN,GS}.
Observe that steps 2 and 3 from Section \ref{proof.that.zero.is.okay} 
can be carried
over easily to the infinite space setting. 
In words, once we have, for any given normal direction, 
a ``good'' hyperplane $T_\lambda$, we can push it along until we
hit a symmetry center.
Here ``good'' means that $\lambda$ lies in the obvious
analogue of (4.3).
Unfortunately we have not been able to
show that there are any good hyperplanes to start at all. 
Finding such a hyperplane would amount to proving that for
$\lambda$ sufficiently negative we have
$u(x^\lambda)-u(x) > 0$ on $\Sigma_\lambda$.
This problem is of course trivial for $u$ on a finite
ball (see step 1 in Section 4.2), and it can be
solved in the classical case by considering the
differential inequality and boundary values one
obtains for $w(x) = u(x^\lambda)-u(x)$ on $\Sigma_\lambda$.
However due to the non-locality of $\Delta_\alpha$
the question becomes more intricate in our case:
Consider without loss of generality a hyperplane $T_0$ 
through the origin and a $w$ that has 
a reflection anti-symmetry $w(x) = -w(x^0)$ and solves
$\Delta_\alpha w(x) + \psi(x) w(x) = 0$, where
$\psi$ is a positive function whose size
and decay behavior can be controlled, given
reasonable decay assumptions on the original $u$.
If we could show that then $w$ does not change sign 
in $\Sigma_0$, we would have found a ``good'' hyperplane.

It is also tempting to conjecture that Theorem \ref{u.is.symmetric} holds
if $\Delta_\alpha$ is replaced by the generator of a 
L\'{e}vy process with rotationally symmetric increment
distribution.
\bigskip

\noindent {\bf Acknowledgement} M. Birkner and A. Wakolbinger would like to thank  
Centro de Investigaci\'{o}n en Matem\'{a}ticas, Guanajuato, Mexico, 
and Erwin Schr\"{o}dinger International 
Institute for Mathematical Physics, Vienna, Austria for their 
kind hospitality. J.A. L\'{o}pez-Mimbela appreciates the kind hospitality of Frankfurt University during his visit in summer 2002, and thanks CONACyT (Mexico) and DAAD (Germany) for partial support.

\end{document}